\documentclass[brochure,english,10pt]{bourbaki}% Pour un expos� en anglais
\usepackage{amssymb,amsfonts,amsmath,footnote,graphicx}
\usepackage[francais]{babel}
\date{Octobre 2012}
\bbkannee{65\`eme ann\'ee, 2012-2013}
\bbknumero{1062}
\title{THE BAUM-CONNES CONJECTURE WITH COEFFICIENTS\\ FOR WORD-HYPERBOLIC GROUPS}
\subtitle{after Vincent Lafforgue}
\author{Michael PUSCHNIGG}
\address{IML, UMR 6206 du CNRS,\\
Universit\'e d'Aix-Marseille\\
13009 Marseille -- France}
\email{michael.puschnigg@\\univ-amu.fr}

\begin{document}
\maketitle

\noindent{\bf INTRODUCTION}

\bigskip
In a recent breakthrough V.~Lafforgue verified the Baum-Connes conjecture with coefficients for all word-hyperbolic groups \cite{41}. This is spectacular progress since it provides the first examples of groups with Kazhdan's property (T) satisfying the conjecture with coefficients \footnote[1] {A proof for the Property (T) groups $Sp(n,1)$ has been announced earlier by P.~Julg in \cite{28}.}. Lafforgue's proof is elementary, but of impressive complexity.\\
\\
In fact, the Baum-Connes conjecture with coefficients is known to be false in general. The first counterexamples were obtained by N.~Higson, V.~Lafforgue, and G.~Skandalis \cite{24} for certain classes of Gromov's random groups \cite{20}. (Note that Gromov's groups are nothing but inductive limits of word-hyperbolic groups !)\\
\\
Already in the early eighties, A.~Connes emphasized that Kazhdan's property (T), which means that the trivial representation of a locally compact group is separated from all other unitary representations, might be a serious obstruction to the Baum-Connes conjecture. The only previously known approach, due to Kasparov \cite{32}, demands the construction of a homotopy among unitary representations between the regular and the trivial representation, which cannot exist for non-compact groups with Property (T). This led to a search for such homotopies among larger classes of representations \cite{26},\cite{35},\cite{41}. V.~Lafforgue \cite{38} introduces the notion of group representations of weak exponential growth. He shows that the trivial representation is not isolated among such representations for hyperbolic groups which opens the way to his proof of the Baum-Connes conjecture with coefficients. For higher rank groups and lattices however, a corresponding version of Property (T) continues to hold \cite{38},\cite{39}. This leads to interesting applications in graph theory and rigidity theory \cite{39} and makes it hard to believe that the Baum-Connes conjecture (at least in the case with coefficients) might be proved for higher rank lattices by the established methods \cite{40}.\\
\\
In section 1, we review index theory and formulate the Baum-Connes conjecture as 
a deep and far reaching generalization of the Atiyah-Singer index theorem. The tools which are used to approach the conjecture are presented in section 2: Kasparov's bivariant $K$-theory \cite{31},\cite{32}, and his construction of ``$\gamma$-elements". Section 3 collects the present knowledge about the Baum-Connes conjecture. In particular, we explain the counterexamples of Higson, Lafforgue, and Skandalis. Section 4 deals with Lafforgue's work on generalizations of Kazhdan's property (T). We discuss his results on his Strengthened Property (T) for higher rank groups and lattices and give an account of their proofs. The applications of his work in graph theory and rigidity theory are mentioned as well. In Section 5 we finally outline V.~Lafforgue's proof of the Baum-Connes conjecture with coefficients for word-hyperbolic groups.\\
\\
{\bf Acknowledgements} \\
I thank Vincent Lafforgue very heartly for his help and advice during the preparation of this manuscript.
It is a pleasure to thank Nigel Higson, Georges Skandalis, and Guoliang Yu for their explanations and constructive remarks.

\section{The Baum-Connes Conjecture}

\subsection{Index theory}

Consider a linear elliptic differential operator $D$ on a smooth compact manifold $M$.
Its {\bf analytical index} is defined as 
$$
Ind_a(D)\,=\,dim\,(Ker D)\,-\,dim(Coker D)\,\in\,{\mathbb Z}.
\eqno(1.1)
$$
The analytical index is invariant under perturbations of the elliptic operator and turns out to be calculable by topological means. In fact, it only depends on the class 
$$
[\sigma_{pr}(D)]\in K^0(T^*M)
\eqno(1.2)
$$ of the principal symbol of $D$. Here $T^*M$ is the total space of the cotangent bundle of $M$ and $K^*$ denotes (compactly supported) topological $K$-theory \cite{3}. (The latter $K$-group can actually be identified with the set of homotopy classes of pseudo-elliptic symbols.) The topological $K$-theory of Atiyah-Hirzebruch is a generalized oriented cohomology theory in the sense of algebraic topology. $K$-oriented manifolds, for example the total space of the cotangent bundle $T^*M$ of a compact manifold $M$, therefore satisfy a $K$-theoretic version of Poincar\'e duality. The image of the symbol class under 
$$
K^0(T^*M)\overset{PD}{\longrightarrow} K_0(M)\overset{p_*}{\longrightarrow} K_0(pt)={\mathbb Z},
\eqno(1.3)
$$ 
$p:M\to pt$ the constant map, is called the {\bf topological index} $Ind_t(D)$ of $D$. 

Suppose now that a compact Lie group $H$ acts smoothly on 
$M$, leaving $D$ invariant.
Then kernel and cokernel of $D$ become finite-dimensional $H$-modules and one may define the equivariant analytical index of $D$ 
$$
Ind_a(D)\,=\,[Ker D]\,-\,[Coker D]\,\in\,R(H),
\eqno(1.4)
$$
as element of the representation ring $R(H)$. The equivariant topological index can be defined in a similar way as before as an element of the equivariant $K$-homology group $K^H_0(pt)$ of a point. There is a tautological isomorphism
$$
\mu: K_0^H(pt)\overset{\simeq}{\longrightarrow} R(H)
\eqno(1.5)
$$
which allows to view both equivariant indices as virtual finite dimensional representations of $H$.
The {\bf Atiyah-Singer Index Theorem} reads then as follows:
\begin{theo} \cite{3}, (6.7). The analytical index and the topological index coincide as homomorphisms  $K_{H}(T^*M) \to R(H)$.
\end{theo}

\subsection{Higher index theory}
Kasparov \cite{32} and Baum-Connes \cite{7},\cite{8} claim that a similar index theorem holds in the following much more general setting:
\begin{itemize}
\item $G$ is an arbitrary locally compact group,
\item $M$ is a smooth manifold equipped with a proper and cocompact $G$-action,
\item $D$ is a $G$-invariant linear elliptic differential operator on $M$.  
\end{itemize}
Note that the condition on the action of $G$ implies that $M$ is non-compact if $G$ is. In particular, $D$ cannot be Fredholm in any naive sense for non-compact $G$. Thus completely new ideas are needed to give a meaning to an ``analytical index".\\

Assume that the locally compact group $G$ acts smoothly and properly on the manifold $M$. Then there exists a $G$-invariant smooth positive measure $dvol$ on $M$. The corresponding Sobolev spaces become $G$-Hilbert spaces, which appear as subrepresentations of a (countable) multiple of the (left)-regular representation on $L^2(G)$. 
\begin{defi}
The {\bf reduced group $C^*$-algebra} of a locally compact group $G$ is the closure in operator norm of the image of the group Banach algebra $L^1(G)$ under the (left)-regular representation:
$$
C^*_r(G)\,=\,\overline{\pi_{reg}(L^1(G))}\,\subset\,\mathcal{L}(L^2(G)).
\eqno(1.6)
$$
\end{defi}
Let $D$ be a $G$-invariant linear elliptic differential operator on $M$. If the $G$-action on $M$ is proper and in addition cocompact one may define an 
{\bf equivariant analytical index} 
$$
Ind^G_a(D)\,=\,``[Ker D]\,-\,[Coker  D]"\,\in\,K_0(C^*_r(G)).
\eqno(1.7)
$$
of $D$. If the kernel and the cokernel of $D$ happen to be finitely generated and projective as modules over $C^*_r(G)$, then the equivariant analytical index of $D$ coincides with their formal difference. As in the classical case the equivariant analytical $C^*$-index is of topological nature and depends only on the symbol class 
$[\sigma_{pr}(D)]\in K_G^0(T^*M).$
The $G$-equivariant topological $K$-theory for {\bf proper} $G$-spaces \cite{47} is very similar to the equivariant $K$-theory with respect to a compact Lie group \cite{3}. In particular, one may define the {\bf equivariant topological index} $Ind^G_t(D)$ of $D$ as the image of the symbol class under
$$
K_G^0(T^*M)\overset{PD}{\longrightarrow} K^G_0(M)\overset{\varphi_*}{\longrightarrow} K^G_0(\underline{E}G),
\eqno(1.8)
$$
where $PD$ denotes K-theoretic Poincar\'e duality and $\varphi:M\to\underline{E}G$ is the equivariant classifying map to a universal proper $G$-space $\underline{E}G$ \cite{8} (such a space always exists and is unique up to equivariant homotopy equivalence). There is a canonical assembly map \cite{7},\cite{8}
$$
\mu: \, K^G_*(\underline{E}G) \, \longrightarrow \, K_*(C^*_r(G)),
\eqno(1.9)
$$
which generalizes (1.5). The corresponding index theorem is 

\begin{theo} \cite{8},\cite{30} Let $G$ be a locally compact group and let $D$ be a $G$-invariant linear elliptic differential operator on the proper, cocompact $G$-manifold $M$. Then
$$
\begin{array}{ccc}
\mu(Ind^G_t(D)) & = & Ind^G_a(D).
\end{array}
\eqno(1.10)
$$
\end{theo}
Every class in $K^G_0(\underline{E}G)$ can be represented by an equivariant topological index, so that the index theorem characterizes the assembly homomorphism $\mu$ as the unique map sending topological to analytical indices.\\
\\
Baum and Connes conjecture that the assembly map provides the link, which allows a purely geometric description of the $K$-theory of reduced group $C^*$-algebras.

\begin{conj}{\bf(Baum-Connes Conjecture)}$\bf{(BC)}$ \cite{8}, (3.15)\\
Let $G$ be a second countable, locally compact group. Then the assembly map 
$$
\begin{array}{cccc}
\mu: & K^G_*(\underline{E}G) & \longrightarrow & K_*(C^*_r(G))
\end{array}
\eqno(1.11)
$$
is an isomorphism of abelian groups.
\end{conj}

\subsection{The conjecture with coefficients}
Baum, Connes, and Higson formulate also a much more general twisted version of conjecture 1.4 \cite{8}. If $D:{\mathcal E}_0\to{\mathcal E}_1$ is a $G$-invariant elliptic operator over the proper and cocompact $G$-manifold $M$, as considered before, then the topological vector spaces ${\mathcal E_0},\,{\mathcal E_1}$ are simultaneously modules over $G$ and the $C^*$-algebra $C_0(M)$ of continuous functions on $M$ vanishing at infinity. One assumes now in addition that

\vspace{3mm} 
\begin{itemize}
\item $\mathcal E_0$ and $\mathcal E_1$ are (right)-modules over an auxiliary $G-C^*$-algebra $A$,
\item The $A$-action on ${\mathcal E}_0,{\mathcal E}_1$ commutes with $D$ and the action of $C_0(M)$,
\item The module multiplications $C_0(M)\otimes{\mathcal E}\to{\mathcal E}$ and ${\mathcal E}\otimes A \to{\mathcal E}$ are $G$-equivariant.
\end{itemize}
\vspace{3mm} 

These conditions imply that the kernel and the cokernel of $D$ are simultaneously $G$-modules and $A$-modules, i.e. they are modules over the following $C^*$-algebra.
\begin{defi}
The {\bf reduced crossed product} of a locally compact group $G$ acting on a $C^*$-algebra $A$ is the closure (in operator norm) of the image of the twisted group Banach algebra $L^1(G,A)$ under the (left)-regular representation:
$$
C^*_r(G,A)\,=\,\overline{\pi_{reg}(L^1(G,A))}\,\subset\,\mathcal{L}(L^2(G,\mathcal{H})),
\eqno(1.12)
$$
$$
(f*\xi)(g)\,=\,\underset{g'g''=g}{\int}\,\pi(g^{-1}\cdot f(g'))\xi(g'')d\mu,
\eqno(1.13)
$$
$$
\,\forall f\in L^1(G,A),\,\,\,\forall \xi \in L^2(G,\mathcal{H}),
$$
where $\pi:\,A\to\mathcal{L}(\mathcal{H})$ is any faithful representation. (The algebra $C^*_r(G,A)$
is independent of the choice of $\pi$.)
\end{defi}
As before, one may define a twisted analytical index 
$$
Ind_a^{G,A}(D)\in K_0(C^*_r(G,A)),
\eqno(1.14)
$$ 
and a twisted topological index 
$$
Ind_t^{G,A}(D)\in K^G_*(\underline{E}G,A).
\eqno(1.15)
$$ 
Here the groups $K^G_*(-,A)$ denote a twisted form of topological $K$-homology for proper $G$-spaces.
Again, there is a corresponding twisted assembly map, which leads to an index theorem with coefficients.

\begin{exem}
If $G=1$ and $A=C(X)$, $X$ a compact Hausdorff space, then $K_*(C^*_r(G,A))\,\simeq\,K^G_*(\underline{E}G,A)\,\simeq\,K^*(X)$ and the previous index theorem equals the index theorem of Atiyah-Singer \cite{4} for families of elliptic operators parametrized by $X$.
\end{exem}

Baum, Connes, and Higson conjecture that the twisted assembly map allows a geometric description of the $K$-theory of reduced crossed product $C^*$-algebras.
\begin{conj} {\bf (Baum-Connes Conjecture with Coefficients)}${\bf (BC_{Coeff})}$,\\ \cite{8}, (6.9).
Let $G$ be a second countable locally compact group and let $A$ be a separable $G$-$C^*$-algebra. Then the assembly map 
$$
\begin{array}{cccc}
\mu_{(G,A)}: & K^G_*(\underline{E}G;A) & \longrightarrow & K_*(C^*_r(G,A))
\end{array}
\eqno(1.16)
$$
from the topological $K$-homology with coefficients in $A$ of a universal 
proper $G$-space $\underline{E}G$ to the $K$-theory of the reduced crossed product 
$C^*$-algebra of $(G,A)$ is an isomorphism of abelian groups.
\end{conj}

\begin{rema} 
For $A={\mathbb C}$ this is just the Baum-Connes conjecture for $G$.
\end{rema}

\begin{rema} 
If ${\bf BC_{coeff}}$ holds for a given group $G$, then it holds for all its closed subgroups $H$.
 More specifically, ${\bf BC_{coeff}}$ for $G$ and $A=C_0(G/H)$ implies ${\bf BC}$ for $H$.
\end{rema}

\section{How to prove the Conjecture}
Classical index theory was not only the point of departure for the developments that led to the Baum-Connes conjecture. Up to now, all attempts to prove it were inspired by Atiyah's index theoretic proof of the Bott periodicity theorem \cite{2}, and rely essentially on Kasparov's bivariant $K$-theory. 

\subsection{Kasparov's bivariant K-theory}

Kasparov's bivariant $K$-theory \cite{31},\cite{32},\cite{16} provides the correct framework and the most advanced technology for the study of (higher) index theory and the $K$-theory of operator algebras. It associates to a pair $(A,B)$ of $C^*$-algebras a ${\mathbb Z}/2{\mathbb Z}$-graded abelian group $KK^*(A,B)$, which is contravariant in $A$ and covariant in $B$. There is a natural isomorphism
$$
\begin{array}{ccc}
K_*(A) & \overset{\simeq}{\longrightarrow} KK^*({\mathbb C},A).
\end{array}
\eqno(2.1)
$$
The bivariant $K$-functor is in both variables
\begin{itemize}
\item {\bf stable}, i.e. it turns the inclusion 
$$
A\hookrightarrow\underset{n\to\infty}{\lim}M_n(A)\simeq A\otimes_{C^*}{\mathcal K}({\mathcal H})
\eqno(2.2)
$$ 
into an isomorphism, and
\item {\bf split exact}, i.e. it maps splitting extensions of $C^*$-algebras into split exact sequences of abelian groups.
\end{itemize}
The key property of Kasparov theory is the existence of a natural associative product  
$$
\begin{array}{cccc}
KK^*(A,B)\otimes KK^*(B,C) & \longrightarrow & KK^*(A,C), 
\end{array}
\eqno(2.3)
$$
making the groups $KK^*(A,A)$ into unital and associative graded algebras. \\

Contrary to ordinary operator $K$-theory, bivariant $K$-theory can be characterized by a simple axiom. 
The Kasparov product allows to define an additive category $\mathcal{KK}$ with (separable) $C^*$-algebras as objects and the even bivariant $K$-groups as morphisms: 
$$
Ob_{\mathcal {KK}}\,=\,C^*-Alg,\,\,\,Mor_{\mathcal {KK}}(A,B)\,=\,KK^0(A,B).
\eqno(2.4)
$$ 
\begin{theo}(Cuntz, Higson)\cite{14},\cite{21}.
Every stable and split exact functor from the category of $C^*$-algebras to an additive category 
factors uniquely through $\mathcal{KK}$.
\end{theo}
In particular, there is a natural transformation
$$
\begin{array}{ccc}
KK^*(A,B) & \longrightarrow & Hom^*(K_*A,\,K_*B).
\end{array}
\eqno(2.5)
$$
For a given locally compact group $G$, there exists an equivariant bivariant $K$-theory $KK_G$ on the category of separable $G$-$C^*$-algebras \cite{32}, which is characterized by a similar universal property \cite{45},\cite{51}. The universal property implies the existence of natural
``Descent" transformations
$$
\begin{array}{ccc}
KK^*_G(A,B) & \longrightarrow & Hom^*(K_*(C^*_r(G,A\otimes C)),\,K_*(C^*_r(G,B\otimes C))),
\end{array}
\eqno(2.6)
$$
which are compatible with the Kasparov product. (Here $C$ is an auxiliary coefficient $C^*$-algebra and the symbol $\otimes$ denotes either the maximal or the minimal $C^*$-tensor product. For commutative $C$, the only case we need, both tensor products coincide.)
\\ 
To apply the theory, one needs an explicit description of bivariant $K$-groups as homotopy classes of $K$-cycles and means to calculate their Kasparov products.
We give such a description in the case $A=B={\mathbb C}$.

\begin{defi} 
(Kasparov \cite{32}) Let $G$ be a (second countable) locally compact group. Then the ring $KK^0_G({\mathbb C},{\mathbb C})$ of Fredholm-representations of $G$ is given by the set of homotopy classes of triples
$$
{\mathcal E}\,=\,({\mathcal H}^\pm,\,\rho^\pm,F),
\eqno(2.7)
$$
where ${\mathcal H}^\pm$ is a ${\mathbb Z}/2{\mathbb Z}$-graded (separable) Hilbert space, equipped 
with an even unitary representation $\rho^\pm$ of $G$, and $F:{\mathcal H}^\pm\to{\mathcal H}^\mp$ is an odd, bounded linear operator such that
$$
\begin{array}{ccc}
F^2\,-\,Id\,\in\,{\mathcal K}({\mathcal H}^\pm) & \text{and} & g\mapsto [F,\rho^\pm(g)]\,\in\,C(G,{\mathcal K}({\mathcal H}^\pm)).
\end{array}
\eqno(2.8)
$$
Here ${\mathcal K}({\mathcal H}^\pm)$ denotes the algebra of compact operators on ${\mathcal H}^\pm$.
\end{defi}

If one writes $F\,=\,\begin{pmatrix} 0 & u \\ v & 0 \end{pmatrix}$, then the conditions (2.8) state that 
$u$ and $v$ are almost equivariant Fredholm operators, which are inverse to each other modulo compact operators. 

If $G$ is compact, then the Fredholm representation ring coincides with the ordinary representation ring. For $G$ abelian, $KK^0_G({\mathbb C},{\mathbb C})$ is 
canonically isomorphic to the topological (Steenrod)-$K$-homology of the dual group $\widehat{G}$, viewed as locally compact topological space.

\subsection{The $\gamma$-element}

All attempts to prove the Baum-Connes conjecture rely up to now on Kasparov's ``Dirac-Dual Dirac" method \cite{32}, which can be viewed as nonlinear version of Atiyah's proof of equivariant Bott-periodicity \cite{2}. 
Suppose for simplicity that there exists a $G$-manifold $M$, which serves as a model for the universal proper $G$-space $\underline{E}G$. Then there exists a canonical class
$$
\alpha\,\in\,KK^G(C_0(T^*M),{\mathbb C}),
\eqno(2.9)
$$
which induces the Baum-Connes map under descent (modulo Poincar\'e-duality). So the assembly map with coefficients factorizes as
$$
\mu_{G,A}:K_*^G(\underline{E}G,A)\,=\,K_*^G(M,A)\,\overset{PD}{\simeq}
K_*(C_r^*(G,C_0(T^*M,A)))\overset{\alpha_*}{\longrightarrow}K_*(C_r^*(G,A)))
\eqno(2.10)
$$
for any $G$-$C^*$-algebra $A$. The key idea is to show that the class $\alpha$ is invertible with respect to the Kasparov product. The Baum-Connes conjecture with coefficients follows then simply by descent. 
In full generality Kasparov's approach to the Baum-Connes conjecture can be summarized as follows:

\begin{theo} ``Dirac-Dual Dirac" Method) \cite{32},\cite{23}
Let $G$ be a locally compact group. Suppose that there exists a locally compact proper $G$-space $X$ and elements $\alpha\in KK^n_G(C_0(X),{\mathbb C})$ and $\beta\in KK^n({\mathbb C}, C_0(X)),\,$ $n$=dim($X$), such that 
$$
\gamma=\beta\otimes\alpha\in KK_G^0({\mathbb C},{\mathbb C})
\eqno(2.11)
$$ 
satisfies
$res^G_H(\gamma)=1\in KK_H^0({\mathbb C},{\mathbb C})$ for every compact subgroup $H$ of $G$. 
Then the Baum-Connes assembly map (with coefficients) for $G$ is split injective. If 
moreover
$$
\gamma=1\in KK_G({\mathbb C},{\mathbb C}),
\eqno(2.12)
$$
or if at least the image of $\gamma$ under descent (2.6) equals the identity,
then the Baum-Connes conjecture (with coefficients) holds for $G$.
\end{theo}
The $\gamma$-element of the previous theorem is unique if it exists \cite{52}. 
\section{Status of the Conjecture}

The Baum-Connes map provides a link between a rather well understood geometric object,
the equivariant $K$-homology of a certain classifying space of a group, and a quite mysterious analytic object, the $K$-theory of its reduced group $C^*$-algebra. The Baum-Connes conjecture appears therefore as quite deep and surprising. It has two aspects: the injectivity of the assembly map (1.11), which has important implications in geometry and topology, and its surjectivity, which proved to be a much more elusive problem. 

The injectivity of the Baum-Connes map with coefficients is known for all connected groups and all groups acting properly and isometrically on a $CAT(0)$-space. Kasparov and Yu \cite{34} recently showed its injectivity for the very huge class of discrete groups, which (viewed as metric spaces with respect to a word metric) admit a uniform coarse imbedding (see (3.10)) into a Banach space $B$ with the following property: there exist an increasing sequence of finite dimensional subspaces of $B$ with dense union, a similar sequence of subspaces of a Hilbert space, and a uniformly continuous family of degree one maps between the corresponding unit spheres of the two families of subspaces. 
\\
A possible source for counterexamples is the seemingly quite different functorial behavior 
of source and target of the Baum-Connes map. Whereas the left hand side of (1.11) is functorial 
under continuous group homomorphisms, this is not at all obvious for the right hand side. 
The reduced group $C^*$-algebra $C^*_r(G)$ is functorial under proper, but not under arbitrary group homomorphisms. For example, the reduced $C^*$-algebra of a non-abelian free group is simple \cite{48}, i.e. has no nontrivial quotients ! It is also easy to see, that the trivial homomorphism $G\to 1$ gives rise to a homomorphism of reduced group $C^*$-algebras iff $G$ is amenable.
The Baum-Connes conjecture claims that the $K$-groups $K_*(C^*_r(G))$ should nevertheless be functorial under arbitrary group homomorphisms, which is quite surprising.\\

At this point one might be tempted to replace the reduced group $C^*$-algebra $C^*_r(G)$ by the maximal group $C^*$-algebra $C^*_{max}(G)$. The representations of the latter correspond to arbitrary unitary representations of $G$ and $C^*_{max}(G)$ is therefore fully functorial in $G$. Examples (see section 4.5) show however that the corresponding assembly map 
$$
\mu_{max}:\,K_*^G(\underline{E}G)\,\longrightarrow\,K_*(C^*_{max}(G))
\eqno(3.1)
$$ 
is far from being an isomorphism in general. 
\\ 
For the conjecture with coefficients, one may study in addition the functoriality with respect to the coefficients of source and target of the Baum-Connes map. This time, the different behavior of both sides leads to the counterexamples to ${\bf BC_{coeff}}$ found by Higson, Lafforgue and Skandalis \cite{24}. We will present them at the end of this section. \\
On the other hand Bost has defined an assembly map
$$
\mu_{L^1}:\,K_*^G(\underline{E}G,A)\,\longrightarrow\,K_*(L^1(G,A))
\eqno(3.2)
$$ 
and conjectures that it is always an isomorphism. This is true for a quite large class of groups \cite{35}. In addition, the counterexamples of \cite{24} do not apply to (3.2).

\subsection{Lie groups and algebraic groups over local fields}
Let $G$ be a connected Lie group and let $H\subset G$ be a maximal compact subgroup. The homogeneous space $G/H$ may serve as a model of the universal proper $G$-space $\underline{E}G$. If $G/H$ carries a $G$-invariant $Spin^c$-structure, the Baum-Connes conjecture equals

\begin{conj} (Connes-Kasparov Conjecture) \cite{8}
Let $i=dim(G/H)\,mod\,2$. Then the map 
$$
\widetilde{\mu}:\,R(H)\,\longrightarrow\,K_i(C^*_r(G)),
\eqno(3.3)
$$
which associates to a virtual representation $[V]\in R(H)$ the $G$-index of the twisted Dirac-operator $\partial_V$ on $G/H$, is an isomorphism of abelian groups. Moreover $K_{i+1}(C^*_r(G))=0$.

\end{conj}
Conjecture 3.1 was proved for linear real reductive groups by A.~Wassermann \cite{53} in 1982. He used many deep results in the representation theory of semisimple Lie groups. 
In his thesis \cite{35}, V.~Lafforgue used geometric methods and employed the existence of a $\gamma$-element to establish ${\bf (BC)}$ for real reductive groups as well as for reductive algebraic groups over local fields. This work was presented at S\'eminaire Bourbaki by G.~Skandalis \cite{50}.

\subsection{Amenable and connected groups}

Following Gromov \cite{19}, a locally compact group is called a-T-menable if it admits a proper, affine, isometric action on a Hilbert space. This is in some sense complementary to Kazhdan's property $(T)$, discussed in the next section. The class of a-T-menable groups contains all amenable groups and all closed subgroups of real and complex Lorentz groups. A proper action of such a group $G$ on a Hilbert space is universal in the sense that the affine Hilbert space may serve as a model for $\underline{E}G$. Higson and Kasparov \cite{23} view an affine Hilbert space as the limit of its finite-dimensional affine subspaces, and use the ``Dirac" and ``Dual Dirac" elements 2.3 on these subspaces to construct a $\gamma$-element $\gamma_G$ for every a-T-menable group $G$. They show that $\gamma=1\in KK^G({\mathbb C},{\mathbb C}),$ and deduce that ${\bf BC_{Coeff}}$ holds for all a-T-menable groups. See the talk of P.~Julg at S\'eminaire Bourbaki \cite{27} for a detailed account to their work. Combining the results of Lafforgue and Higson-Kasparov, Chabert, Echterhoff and Nest \cite{12} succeeded finally in verifying ${\bf BC}$ for all locally compact, connected groups.

\subsection{Discrete groups}
Let $G=\Gamma$ be a countable discrete group. We suppose for simplicity that $\Gamma$ is torsion free. Any contractible, proper and free $\Gamma$-space $E\Gamma$ may serve as a model for $\underline{E}\Gamma$ and the Baum-Connes conjecture equals 

\begin{conj} \cite{7} Let $\Gamma$ be a torsion free, countable discrete group and let $B\Gamma$ be a classifying space for principle-$\Gamma$-bundles. Then the assembly map
$$
\mu:\,K^{top}_*(B\Gamma)\,\longrightarrow\,K_*(C^*_r(\Gamma))
\eqno(3.4)
$$
is an isomorphism of abelian groups.
\end{conj}

The most important progress up to now was achieved by V.~Lafforgue \cite{35},\cite{36}, who
established ${\bf BC}$ for word-hyperbolic groups in the sense of Gromov and for uniform lattices in the higher rank groups $SL_3(K)$, $K$ a local field. He and P.~Julg \cite{28} were the first who overcame  the barrier of Kazhdan's Property $(T)$ (which holds for generic hyperbolic groups and all higher rank lattices). For both classes of groups there exists a $\gamma$-element, but it cannot be equal to 1 in the presence of property $(T)$. Nevertheless $\gamma$ acts as the identity on $K_*(C^*_r(\Gamma))$ which already implies ${\bf BC}$. See also Skandalis' report at S\'eminaire Bourbaki \cite{50}.

\subsection{The conjecture with coefficients}

The Baum-Connes conjecture with coefficients was previously known only for a-$T$-menable groups by \cite{23}, and for hyperbolic groups and commutative (!) coefficients by the work of 
Lafforgue \cite{37}. The first proof of ${\bf BC_{Coeff}}$ for a class of groups with Property (T), the Lie groups $Sp(n,1)$, is due to Julg and sketched in \cite{28}. 
The spectacular recent breakthrough, which is the main topic of this expos\'e, is again due to Vincent Lafforgue:

\begin{theo} (Lafforgue) \cite{41}
The Baum-Connes conjecture with coefficients holds for all locally compact groups acting 
properly and isometrically on a weakly geodesic and locally uniformly finite 
hyperbolic metric space. In particular, it holds for all word-hyperbolic groups.
\end{theo}

Contrary to the Baum-Connes conjecture, which is open at the moment, the Baum-Connes conjecture with coefficients is known to be false in general. 

\subsection{A counterexample}

In recent years Gromov's spectacular theory of ``Random Groups" \cite{20},\cite{17} has been used to produce various counterexamples to open questions in geometric group theory and operator algebras. One instance is the following counterexample to the Baum-Connes conjecture with coefficients, which is due to Higson, Lafforgue, and Skandalis \cite{24}. It is based on the possibility of embedding some expander graphs coarsely and uniformly into the Cayley graphs of random groups.
\\
As indicated before, it is the different functorial behavior of source and target of the Baum-Connes 
assembly map $\mu_{(G,A)}$, which leads to the desired counterexamples. The main idea is the following.
Let $\Gamma$ be a discrete group. Suppose that there exists an extension
$$
0\,\to\,I\,\to\,A\,\to\,B\,\to\,0
\eqno(3.5)
$$
of $\Gamma$-$C^*$-algebras, ($I\subset A$ an ideal and $B\simeq A/I$), such that the upper line in the commutative diagram 
$$
\begin{array}{ccccc}
K_*(C^*_r(\Gamma,I)) & \to & K_*(C^*_r(\Gamma,A)) & \to & K_*(C^*_r(\Gamma,B)) \\

\uparrow & & \uparrow & & \uparrow\\

 K^{\Gamma}_*(\underline{E}\Gamma,I) & \to & K^{\Gamma}_*(\underline{E}\Gamma,A) & 
\to & K^{\Gamma}_*(\underline{E}\Gamma,B) \\
\end{array}
\eqno(3.6)
$$ 
is not exact in the middle. As the lower line is always exact in the middle, one deduces that the vertical arrows,
given by the corresponding Baum-Connes assembly maps, cannot all be isomorphisms, as the Baum-Connes conjecture with
coefficients predicts.\\
The key point is therefore to find a projector $p\in C^*_r(\Gamma,A)$, whose class in $K$-theory is not in the image of $K_*(C^*(\Gamma,I))$, and which maps to 0 in $C^*_r(\Gamma,B)$.

Recall that the Laplace operator $\Delta:\ell^2(\mathcal {G})\to\ell^2(\mathcal {G})$ on a graph $\mathcal G$ of bounded valency is the positive bounded linear operator given by 
$$
\begin{array}{ccc}
\Delta f(x) & = &  \underset{d(x,y)=1}{\sum} \, \left( f(x)-f(y) \right),
\end{array}
\eqno(3.7)
$$
where the sum runs over the set of all vertices $y$ adjacent to $x$.
 If $\mathcal G$ is finite, then the kernel of $\Delta$ coincides with the space of locally constant functions. \\
A sequence $(\mathcal{G}_n),n\in{\mathbb N},$ of finite connected graphs of uniformly bounded valency is an 
{\bf expander} \cite{15},\cite{9}, if their cardinality tends to infinity 
$$
\underset{n\to\infty}{\lim}\,\sharp({\mathcal{G}_n})\,=\,\infty,
\eqno(3.8)
$$
and if their Laplace operators have a uniform spectral
 gap. i.e. $\exists\,\epsilon>0:$
$$
\begin{array}{cccccc}
 Sp(\Delta(\mathcal{G}_n)) & \cap & ]0,\epsilon[ & = & \emptyset, & \forall n\in\mathbb N.
\end{array}
\eqno(3.9)
$$
The Cayley graph ${\mathcal G}(\Gamma,S)$ of a finitely generated group $(\Gamma,S)$ has the group $\Gamma$ itself as set of vertices,
and two vertices $g,h\in\Gamma$ are adjacent iff $g^{-1}h\in S$.

Now, according to Gromov \cite{20},\cite{17},\cite{1}, it is possible to imbed a suitable expander coarsely and uniformly
 into the Cayley graph of some finitely generated group $(\Gamma,S)$. This means 
that there exists a sequence $i_n:\mathcal{G}^0_n\to {\mathcal G}^0(\Gamma,S),\,n\in{\mathbb N},$ of maps of vertex sets, such that 
$$
\begin{array}{ccccccc}
 \rho_0(d_{\mathcal{G}_n}(x,y)) & \leq & d_{{\mathcal G}(\Gamma,S)}(i_n(x),i_n(y)) & \leq & \rho_1(d_{\mathcal{G}_n}(x,y)),& \forall x,y\in\mathcal{G}_n, & \forall n\in\mathbb N,
\end{array}
\eqno(3.10)
$$
for some monotone increasing, unbounded functions $\rho_0,\,\rho_1:{\mathbb R}_+\to{\mathbb R}_+$.
The coarse imbeddings $i_n$ (which we suppose to be injective to simplify notations) 
may be used to ``transport" the Laplace operators of the expander graphs 
to an operator on $\ell^2(\Gamma)$.To be precise let 
$$
\Theta_n:\ell^2(\mathcal{G}_n)\to\ell^2(\Gamma),\,e_x\,\mapsto\,e_{i_n(x)}
\eqno(3.11)
$$
and put 
$$
\Delta'_n\,=\,\Theta_n\Delta(\mathcal{G}_n)\Theta_n^*\,+\,(1-\Theta_n\Theta^*_n)\,\in\,
{\mathcal L}(\ell^2(\Gamma)).
\eqno(3.12)
$$ 
Consider the operator $\Delta'=\underset{n}{\bigoplus}\,\Delta'_n$ on the Hilbert sum 
$\mathcal{H}\,=\,\underset{n}{\bigoplus}\,\ell^2(\Gamma)\,=\,\ell^2(\mathbb N\times\Gamma).$
\\
The reduced crossed product $C^*_r(\Gamma,C_b(\mathbb N, C_0(\Gamma)))$ acts faithfully on $\mathcal{H}$.
 The operator $1-\Delta'$ may be written as a finite (!) sum $\underset{g}{\sum} f_g u_g,\,f_g\in C_b(\mathbb N, C_0(\Gamma))$
 because of (3.10) and the fact that the propagation speed of the Laplace operator on a graph is equal to one.
 In particular $1-\Delta'\in C_c(\Gamma,C_b(\mathbb N, C_0(\Gamma)))\subset
C^*_r(\Gamma, C_b(\mathbb{N},C_0(\Gamma)))$. It is a positive operator which, according to (3.9), has a spectral gap,
 i.e. $Sp(\Delta')\,\cap\, ]0,\epsilon[ \,= \, \emptyset$. The spectral projection $p'\,=\,\underset{n}{\bigoplus}\,p'_n$ onto 
$$
Ker(\Delta')\,=\,\underset{n}{\bigoplus}\,Ker(\Delta'_n)\,=\,\underset{n}{\bigoplus}\,\mathbb{C}
\eqno(3.13)
$$ 
may thus be obtained from $\Delta'$ by continuous functional calculus, so that
$$
p'\in C^*_r(\Gamma,C_b(\mathbb N, C_0(\Gamma))).
\eqno(3.14)
$$ 
This is the projection we are looking for. As an element of the reduced crossed product,
it can be uniquely written as infinite sum $p'\,=\,\underset{g}{\sum}\,f_g\,u_g,\,f_g\in C_b(\mathbb N, C_0(\Gamma)),\,g\in\Gamma.$ 
It follows from (3.12) and (3.8) that
$$
f_g\in C_0(\mathbb N, C_0(\Gamma)),\,\forall g\in\Gamma.
\eqno(3.15)
$$
Consider now the extension
$$
0\,\to\,C_0(\mathbb N, C_0(\Gamma))\,\to\,C_b(\mathbb N, C_0(\Gamma))\,\to\,Q\,\to\,0
\eqno(3.16)
$$
of $\Gamma$-$C^*$-algebras. 
On the one hand, the image of the projection $p'\in C^*_r(\Gamma,C_b(\mathbb N, C_0(\Gamma)))$ in 
$C^*_r(\Gamma,Q)$ is zero by (3.15). On the other hand, its $K$-theory class 
$$
[p']\in K_0(C^*_r(\Gamma,C_b(\mathbb N, C_0(\Gamma))))
\eqno(3.17)
$$ 
does not lie in the image of 
$K_0(C^*_r(\Gamma,C_0(\mathbb N, C_0(\Gamma))))\,=\,\underset{\underset{n}{\to}}{lim}\,K_0(C^*_r(\Gamma,C_0(\Gamma)))$ 
because 
$$
\pi_n([p'])\,=\,[p'_n]\,\neq\,0\,\in\,K_0(C^*_r(\Gamma,C_0(\Gamma)))\,=\,K_0(\mathcal{K}(\ell^2(\Gamma)))\,
\simeq\,\mathbb{Z},\,\,\forall n\in\mathbb{N}.
\eqno(3.18)
$$
In this way Higson, Lafforgue and Skandalis obtain the desired counterexample.

\section{Kazhdan's Property (T) and its generalizations}

The most important classes of groups, for which the Baum-Connes conjecture is unsettled, are 
simple linear groups of split rank $\geq 2$ over local fields, where ${\bf BC_{Coeff}}$ is open, and lattices in such groups where already ${\bf BC}$ is unknown in most cases. These classes are distinguished by their astonishing rigidity properties \cite{42}. They also provide the most prominent examples of groups with Kazhdan's property (T), which plays a key role in rigidity theory, and has important applications in operator algebras, representation theory, and graph theory \cite{15},\cite{9},\cite{42}. In this section we report what is known about V.~Lafforgue's strengthened versions of Property (T) \cite{38}, and outline its applications to graph theory and rigidity theory \cite{39}. Strengthened Property (T) appears also to be very serious obstruction against a possible ``Dirac-Dual Dirac" approach to the Baum-Connes conjecture \cite{40}.

\subsection{Property (T)}
Recall that a locally compact group $G$ has {\bf Kazhdan's Property (T)} if the following equivalent conditions hold:
\begin{itemize}
\item The trivial representation is an isolated point in the unitary dual of $G$.
\item Every unitary representation $\pi$ of $G$ with almost invariant vectors, i.e.
$$
\forall\epsilon>0,\,\forall K\subset G\,\text{compact},\,\exists\,\xi\in\mathcal{H}-\{0\}:\,\parallel \pi(g)\xi-\xi\parallel\leq\epsilon\parallel\xi\parallel,\,\forall g\in K,
$$
contains nonzero fixed vectors.
\item Every continuous isometric affine action of $G$ on a Hilbert space has a fixed point.
\item There exists a projection $p\in C^*(G),$ such that for any unitary representation $(\pi,\mathcal{H})$ of $G$ 
the operator $\pi(p)\in\mathcal{L}(\mathcal{H})$ equals the orthogonal projection onto $\mathcal{H}^{\pi(G)}$.
\end{itemize}
Here $C^*(G)$ denotes the ``full'' group $C^*$-algebra, i.e. the enveloping $C^*$-algebra of the involutive Banach algebra $L^1(G)$.\\
\\
{\bf Examples of Kazhdan groups:}
\begin{itemize}
\item Compact groups, 
\item Simple algebraic groups of split rank at least two over local fields and their lattices.
\item Many hyperbolic groups, for example lattices in the simple Lie groups $Sp(n,1)$,
$n>1,$ or $F_{4,-20}$ of real rank one.
\item Generic, randomly produced hyperbolic groups \cite{20}.
\end{itemize} 
Basic examples of groups without Property (T) are free groups and non-compact amenable or a-T-menable groups.

\subsection{Lafforgue's Strengthened Property (T)}

In recent years various generalizations of Property (T) have been proposed. These deal with larger classes of representations than the unitary ones. A first example is

\begin{defi}(Bader, Furman, Gelander, Monod, \cite{6})\\ 
a) A locally compact group $G$ has ${\bf Property (T)_{uc}}$, if every isometric representation of $G$ on a uniformly convex Banach space with almost invariant vectors has non zero fixed vectors.\\
b) A locally compact group $G$ has ${\bf Property (F)_{uc}}$, if every affine isometric action of $G$
on a uniformly convex Banach space has a fixed point.
\end{defi}
Lafforgue goes one step further and allows not only isometric representations, but representations of weak exponential growth.

\begin{defi} (Lafforgue)
Let $G$ be a locally compact group with a proper, continuous and symmetric length function $\ell:G\to\mathbb{R}_+,$ and let $\lambda > 1$. A continuous representation $\pi$ of $G$ on a Banach space $B$ is of {\bf exponent} ${\bf\lambda}$ (with respect to $\ell$) if
$$
\parallel\pi\parallel_{\lambda}\,=\,\underset{g\in G}{sup}\,\,\,\lambda^{-\ell(g)}\parallel\pi(g)\parallel_{\mathcal{L}(B)}\,<\,\infty.
\eqno(4.1)
$$
\end{defi}
The representations of $G$ of exponent $\lambda$ on a self-dual class of Banach spaces $\mathcal B$ give rise to representations of the corresponding involutive Fr\'echet algebra \cite{36}
$$
\mathcal{C}_\lambda(G,\ell,\mathcal{B}),
\eqno(4.2)
$$
obtained by completion of the convolution algebra of compactly supported continuous functions
on $G$ with respect to the seminorms
$$
\parallel f\parallel_N\,=\,\underset{\underset{\parallel\pi\parallel_\lambda\leq N}{(\pi,B),}}{\sup}\,
\parallel\pi(f)\parallel_{\mathcal{L}(B)},\,N\,\in\,{\mathbb N}.
\eqno(4.3)
$$
The supremum is taken over all representations $(\pi,B)$ of exponent $\lambda>1$ and $\lambda$-norm\\ $\parallel\pi\parallel_\lambda\leq N$ on a space $B\in\mathcal{B}$.

\begin{defi} (Lafforgue),\cite{38},\cite{39}.\\
a) Let $\mathcal B$ be a class of Banach spaces which is closed under taking duals. A locally compact group $G$ has Lafforgue's ${\bf Property \, (T^{Strong}_{\mathcal{B}})}$, if for every proper symmetric length function $\ell$ on $G$, there exists $\lambda>1$ and a selfadjoint idempotent $p_\lambda\in\mathcal{C}_\lambda(G,\ell,\mathcal{B})$ such that 
$$
\begin{array}{ccc}
\pi(p_\lambda)(B) & = & B^{\pi(G)}
\end{array}
\eqno(4.4)
$$
for every representation $(\pi,B)$ of exponent $\lambda$ on a Banach space $B\in\mathcal{B}$.
Such an idempotent is unique and central in $\mathcal{C}_\lambda(G,\ell,\mathcal{B})$.\\
b) It satisfies $\bf Property \, (T)_{Hilb}^{Strong}$ if $\bf (T^{Strong}_{\mathcal{B}})$ holds for 
$\mathcal B\,=\,\{\text{Hilbert spaces}\}.$\\
c)  It possesses $\bf Property \, (T)_{Ban}^{Strong}$ if $\bf (T^{Strong}_{\mathcal{B}})$ holds for every class $\mathcal B$ which is (uniformly) of type $>$ 1. This means that there exist $n\in\mathbb{N}$ and $\epsilon>0$ such that no $n$-dimensional subspace of any $B\in\mathcal{B}$ is $(1+\epsilon)$-isometric to $\ell^1_n$.
\end{defi}
Note that every uniformly convex space is of type $>1$.
The relations between these properties are displayed below.

$$
\begin{array}{ccccc}
{\bf (T)_{Ban}^{Strong}} & \Rightarrow & {\bf (F)_{uc}} & \Rightarrow & {\bf (T)_{uc}} \\
 & & & & \\
\Downarrow & & & & \Downarrow \\
 & & & & \\
{\bf (T)_{Hilb}^{Strong}} & & \Rightarrow & & {\bf (T)}\\
\end{array}
\eqno(4.5)
$$

\subsection{Results}

Concerning the strengthened property (T) one observes a strict dichotomy between groups of ``split rank'' one and ``higher rank'' groups. Despite the fact that they generically satisfy the ordinary Kazhdan property, word-hyperbolic groups are very far from sharing the strengthened versions of property (T).

\begin{theo} (Lafforgue, \cite{38})
Word-hyperbolic groups do not satisfy $\bf (T)_{Hilb}^{Strong}$.
\end{theo}

Lafforgue's proof is closely linked to his work on the Baum-Connes conjecture and will be explained in section 5. The following remarkable result of Yu asserts that hyperbolic groups do not have property 
$\bf (F)_{uc}$ either.

\begin{theo} (Yu,\cite{54})
Every hyperbolic group admits a proper, isometric, affine action on an $\ell^p$-space for $p\in]1,\infty[$ sufficiently large.
\end{theo}

Yu's construction of the desired affine action is related to an explicit description of the $\gamma$-element of a hyperbolic group. \\

For higher rank groups and their lattices however, many (and conjecturally all) strengthened versions of the Kazhdan property hold. 

\begin{theo} (Lafforgue, \cite{38})
A simple real Lie group, whose Lie algebra contains a copy of $\mathfrak{sl}_3$, satisfies $\bf (T)_{Hilb}^{Strong}$. The same holds for its uniform lattices.
\end{theo}

\begin{theo} (Lafforgue, \cite{39})
A simple linear algebraic group over a non-archimedian local field, whose Lie algebra contains a copy of $\mathfrak{sl}_3$, satisfies $\bf (T)_{Ban}^{Strong}$. The same holds for its uniform lattices.
\end{theo}
This result has applications in graph theory.
It is well known that an expanding sequence of graphs (3.8), (3.9) cannot be imbedded uniformly (3.10) into Hilbert space. 
\begin{theo} (Lafforgue, \cite{39})
Let $(\Gamma,S)$ be a uniform lattice in a simple algebraic group over a non-archimedian local field, whose Lie algebra contains a copy of $\mathfrak{sl}_3$. Let $(\Gamma_n),n\in\mathbb{N},$ be a sequence of finite index normal subgroups of $\Gamma$, whose intersection is 1. Then the sequence of (finite) Cayley graphs $(\mathcal{G}(\Gamma/\Gamma_n),\pi(S))$ cannot be imbedded uniformly in any Banach space of type $>1$.
\end{theo}
Recently, Mendel and Naor  \cite{43},\cite{44} used completely different methods 
to construct huge families of expanders which do not admit a uniform embedding into any uniformly convex Banach space.

\subsection{Proofs}

Let $G=SL_3(F)$, $F$ a local field. Let $K$ be a maximal compact subgroup of $G$. Lafforgue's key observation (which generalizes the Howe-Moore property of unitary representations \cite{25}) is that the matrix coefficients of $K$-invariant vectors in representations of sufficiently small exponent tend very quickly (exponentially fast) to a limit at infinity:
$$
\vert\langle\xi,\pi(g)\eta\rangle\,-\,c_{\xi,\eta}\vert\,=\,O(e^{-\mu(\lambda)\ell(g)}),
\eqno(4.6)
$$
where $(\pi,B)$ is a representation of $G$ of exponent $\lambda$, $\xi\in B^*,\,\eta\in B$ are 
$K$-invariant vectors, and $\mu(\lambda)>0$ if $\lambda$ is close to 1. 
(Here $B$ is a Hilbert space in the archimedian case and a Banach space in a class $\mathcal{B}$ of type $>1$ in the non-archimedian case.) 

It is then easy to see that for a fixed, compactly supported positive function of mass one $\chi\in C_c(G)$  the family 
$$
f_g:\,x\mapsto\,\underset{K\times K}{\int}\,\chi(k_1gxk_2)\,dk_1dk_2,\,g\in G,
\eqno(4.7)
$$
of $K$-biinvariant, compactly supported functions on $G$ tends to a limit 
$$
\underset{g\to\infty}{\lim}\,f_g\,=\,p_\lambda\,\in\,\mathcal{C}_\lambda(G,\ell,\mathcal{B})
\eqno(4.8)
$$ 
as $g\in G$ tends to infinity. It follows from a non-spherical version of (4.6) that $p_\lambda$ is the desired ``Kazhdan"-projection. It is selfadjoint as limit of selfadjoint functions. 
This establishes 4.6 and 4.7.

We outline Lafforgue's strategy for proving the decay estimates (4.6) in the non-archimedian case.\\
\\
Let $F$ be a non-archimedian local field with ring of integers $\mathcal O$ and residue field $\mathbb{F}_q$. Let $G=SL(3,F)$ and put $K=SL(3,\mathcal{O})$. There is a Cartan decomposition
$G=K\overline{A}_+K$ with
$$
\overline{A}_+\,=\,\{\,diag(\pi^{-i_1},\pi^{-i_2},\pi^{-i_3}),\,i_1+i_2+i_3=0,\,i_1\geq i_2\geq i_3\},
\eqno(4.9)
$$
where $\pi$ denotes a fixed uniformizing element of $F$. A canonical $K$-biinvariant proper length function on $G$ is given by 
$$
\ell(kak')\,=\,i_1(a)-i_3(a).
\eqno(4.10)
$$
Let $\mathcal B$ be a class of type $>1$ of Banach spaces, which is closed under taking duals. Let $(\pi,B),\,B\in\mathcal{B},$ be a representation of $G$ of exponent $\lambda$ and denote by $(\check{\pi},B^*)$ its contragredient representation. Let $\eta\in B,\xi\in B^*$ be $K$-invariant unit vectors. The corresponding matrix coefficient $g\mapsto\langle\xi,\pi(g)\eta\rangle$ is then determined by its values
$$
c(i_1-i_2,i_2-i_3)\,=\,\langle\xi,\pi(diag(\pi^{-i_1},\pi^{-i_2},\pi^{-i_3}))\eta\rangle,\,i_1+i_2+i_3=0,\,i_1\geq i_2\geq i_3.
\eqno(4.11)
$$
Fix integers $m\geq n\geq 0, m+n \in 3{\mathbb N}$. Lafforgue finds two finite families
 $(\overline{a}_i)_{i\in I},\,(\overline{b}_j)_{j\in J}$ of elements of $G/K$ (considered as points of the affine building), and a matrix
$T\in M_{IJ}({\mathbb C})$ satisfying
\begin{itemize}
\item[a)] $\ell(\overline{a}_i)\,=\,m,\,\,\,\ell(\overline{b}_j)\,=\,n,\,\,\,\forall i\in I,\,\forall j\in J.$\\
\item[b)] $\vert I\vert\,=\,q^{2m},\,\,\,\vert J\vert\,=\,q^{2n}.$\\
\item[c)] The Schur product $\widetilde{T}\in M_{IJ}({\mathbb C}),\,\widetilde{T}_{ij}\,=\,T_{ij\,}\langle\check{\pi}(\overline{a_i})\xi,\,\pi(\overline{b_j})\eta\rangle,$ satisfies
$$
\vert I\vert^{-\frac12}\,\vert J\vert^{-\frac12}\,\underset{i,j}{\sum}\,\widetilde{T}_{ij\,}\,=\,c(m-n+2,n-1)\,-\,c(m-n,n).
\eqno(4.12)
$$
\item[d)] The norm of the operator 
$ T\otimes id_B:\,L^2(J,B)\to L^2(I,B)$ is bounded by the norm of the (normalized) Fourier transform 
$$
\mathcal{F}_A^B:L^2(A,B)\to L^2(\check{A},B),\,\,\,f\mapsto(\chi\mapsto \frac{1}{\sharp A}\,\underset{a}{\sum}
\,\chi(a)f(a))
\eqno(4.13)
$$
on a finite abelian group $A=A_J$ of order $\vert J\vert^{\frac12}$.
\end{itemize}

This allows Lafforgue to bring Fourier analysis on finite abelian groups and the geometry of Banach spaces into play. According to Bourgain \cite{11}, the (normalized) Fourier transform 
satisfies a uniform bound of the type
$$
\parallel\mathcal{F}_A^B\parallel\,=\,O((\sharp A)^{-\alpha})),\,\,\,\alpha\,=\,\alpha({\mathcal B})\,>\,0,
\eqno(4.14)
$$
for every finite abelian group $A$ and every Banach space $B$ in a class $\mathcal{B}$ of type $>1$.
\\
Lafforgue derives thus from (4.12) the estimate
$$
\vert c(m-n+2,n-1)-c(m-n,n)\vert
$$
$$
\leq\,\parallel\widetilde{T}\parallel\,\leq\,\parallel{\mathcal F}^B_{A_J}\parallel\left(
\underset{i\in I}{\max}\parallel\check{\pi}(\overline{a}_{i})\xi\parallel\right)\left(\underset{j\in J}{\max}
\parallel\pi(\overline{b}_{j})\eta\parallel\right)
$$

$$
\leq\,\parallel\pi\parallel^2_\lambda\,(q^{-\alpha})^{n}\,\lambda^{m+n}.
\eqno(4.15)
$$
This, together with the analogous estimate obtained by exchanging the roles of $m$ and $n$, implies for $\lambda>1$ sufficiently close to 1 the exponential decay of differences of matrix coefficients.
Claim (4.6) follows then by a simple Cauchy sequence argument.
\subsection{Relation to the Baum-Connes conjecture}

It was realized very early by A.~Connes, that Kazhdan's Property (T) might be a serious obstruction against the validity of the Baum-Connes conjecture for a noncompact group. At least Kasparov's original ``Dirac-Dual Dirac" method cannot possibly work in the presence of Property (T).\\ 
To see this, recall that the unitary representations of a locally compact group $G$ correspond bijectively to the representations of its full group $C^*$-algebra $C^*(G)$.
In particular, there are epimorphisms $\pi_{reg}:\,C^*(G)\to C^*_r(G)$ and $\pi_{triv}:\,C^*(G)\to\mathbb{C}$ corresponding to the regular and the trivial representation of $G$, respectively.\\
Now Connes argues as follows. \\

For every locally compact group $G$ one may construct an assembly map with values in the $K$-theory of the full group $C^*$-algebra. It fits into the commutative diagram
$$
\begin{array}{cccc}
\mu_{max}: & K_*^{top}(\underline{E}G) & \longrightarrow & K_*(C^*(G)) \\
& \parallel & & \downarrow\,\pi_{reg *} \\
\mu: & K_*^{top}(\underline{E}G) & \longrightarrow & K_*(C^*_r(G)). \\
\end{array}
\eqno(4.16)
$$
Suppose that there exists a $\gamma$-element for $G$, which equals one: $\gamma\,=\,1\,\in KK^G({\mathbb C},\mathbb{C})$. Then both assembly maps have to be isomorphisms and one deduces that $\pi_{reg *}$ is a bijection. If $G$ has Kazhdan's property (T) while not being compact, then the class $[p]\in K_0(C^*(G))$ of the Kazhdan projection is nontrivial because $\pi_{triv}([p])=1\in K_0(\mathbb{C})=\mathbb{Z}$, but maps to zero in $K_0(C^*_r(G))$ because the regular representation of a non compact group has no fixed vectors. Thus $\gamma\neq 1$ for a non compact group with property $(T)$. \\
A beautiful argument of Skandalis \cite{49} shows that for hyperbolic groups with Property (T) even the image of $\gamma$ under descent to bivariant $K$-theory \cite{32} differs from 1: $$j_r(\gamma)\,\neq\,j_r(1)\,=\,1\,\in\,KK(C^*_r(\Gamma),\,C^*_r(\Gamma)).
\eqno(4.17)$$

Nevertheless, it is sometimes possible to show that $\gamma$ maps to the identity under descent (2.6) even for  property (T) groups. The idea, originally due to Julg \cite{26}, is to find enlarged versions $\widetilde{KK}$ of bivariant $K$-theory, which will not have particularly nice properties, but allow to factorize the descent map as 
$$
KK_\Gamma({\mathbb C},{\mathbb C})\,\to\,\widetilde{KK}_\Gamma({\mathbb C},{\mathbb C})\,\to\,
Hom(K_*(C^*_r(\Gamma,A)),K_*(C^*_r(\Gamma,A))),
\eqno(4.18)
$$
and satisfy
$$
[\gamma]\,=\,[1]\,\in\,\widetilde{KK}_\Gamma({\mathbb C},{\mathbb C}).
\eqno(4.19)
$$
In \cite{35}, Lafforgue developed a bivariant $K$-theory for Banach algebras to deal at least with $\bf (BC)$. In the case with coefficients, the absence of $\bf (T)_{Hilb}^{Strong}$ for word-hyperbolic groups (4.4) enables Lafforgue to construct the desired homotopy between $\gamma$ and 1 using bivariant $K$-cycles, whose underlying representations are of small exponential growth \cite{41}.
For general higher rank lattices property $\bf (T)_{Hilb}^{Strong}$ is a very serious obstruction against an implementation of the ``Dirac-Dual Dirac" approach. Lafforgue explains in 
\cite{36},\cite{40} that the only known way to establish (4.18) and (4.19) consists in finding a homotopy between $\gamma$ and $1$ among representations which define bounded Schur multipliers on some isospectral subalgebra of $C^*_r(\Gamma)$. For lattices in $SL_3(F),\,F$ a local field, the existence of such an algebra would contradict (4.12). In fact, it was this circle of ideas which led Lafforgue to the invention of Strengthened Property (T).

\section{Lafforgue's Approach}

In this last section we discuss Lafforgue's proof of the Baum-Connes conjecture with coefficients for word-hyperbolic groups \cite{41}. Recall that a geodesic metric space $(X,d)$ is {\bf hyperbolic}
\cite{18}, if there exists $\delta\geq 0,$ such that for any points $a,b,c\in X$ 
$$
x\,\in\,geod(a,c)\,\Rightarrow\,d(\,x,\,geod(a,b)\cup geod(b,c)\,)\,<\,\delta,
\eqno(5.1)
$$
where for any (nonempty) subsets $A,B\subset X$ 
$$
geod(A,B)\,=\,\{x\in X,\,\exists a\in A,\,\exists b\in B,\,d(a,x)+d(x,b)=d(a,b)\}
\eqno(5.2) 
$$
denotes the union of all geodesic segments joining a point of $A$ and a point of $B$.
A finitely generated group $\Gamma$ is word-hyperbolic \cite{18}, if for one (and therefore every) finite symmetric set of generators $S$ the Cayley-graph $\mathcal{G}(\Gamma,S)$ is a hyperbolic metric space. An important class of word-hyperbolic groups is provided by fundamental groups of compact Riemannian manifolds of strictly negative sectional curvature.\\
\\
There is a distinguished class of models for the universal proper $\Gamma$-space $\underline{E}\Gamma$ of 
a hyperbolic group: the {\bf Rips complexes} \cite{18}. For fixed $R>0,$ the Rips complex $\Delta^R_*(\Gamma,d_S)$ is the simplicial set of finite, oriented subsets of $\Gamma$ of diameter at most $R$:
$$
S'\,\in\,\Delta^R_p(\Gamma,d_S)\,\Leftrightarrow\,S'\subset\Gamma,\,\vert S'\vert=p+1,\,diam(S')\leq R.
\eqno(5.3)
$$
(Here and in the sequel we will use the same notation for a Rips-simplex and its underlying set.) The natural action of $\Gamma$ on $\Delta^R_*(\Gamma,d_S)$ induced by left translation is simplicial and proper. For hyperbolic groups the Rips complex is in addition contractible, provided that $R$ is sufficiently large. It may therefore serve as model for $\underline{E}\Gamma$.
The associated chain complex 
$$
(\mathbb{C}(\Delta^R_*(\Gamma, d_S)),\,\partial),\,\,\,
\partial(g_0,\ldots,g_n)\,=\,\underset{i=0}{\overset{n}{\sum}}\,(-1)^i\,(g_0,\ldots,\widehat{g}_i,
\ldots,g_n)
\eqno(5.4)
$$
is a $\Gamma$-finite and $\Gamma$-free resolution of the constant $\Gamma$-module $\mathbb C$.\\
\\
Various authors \cite{33},\cite{35},\cite{41} have constructed $\gamma$-elements (2.11) for hyperbolic groups. 

\begin{theo}(Kasparov, Skandalis, \cite{33})
For a suitable choice of a hyperbolic distance $d'$ on $\Gamma$, a square zero contracting chain homotopy $h$ (see 5.2) of $({\mathbb C}(\Delta^R_*),\partial)$, and $R,t>>0$ sufficiently large
$$
\left(\ell^2(\Delta^R_*),\,e^{t d'_{x_0}}\,(\partial\,+\,h)\,e^{-t d'_{x_0}}\right)
\eqno(5.5)
$$
defines a bounded $K$-cycle representing $\gamma\in KK^\Gamma({\mathbb C},{\mathbb C})$.
\end{theo}

The $K$-cycle (5.5) is in fact a slightly modified version of the original $\gamma$-element of Kasparov and Skandalis. It is better adapted to Lafforgue's needs \cite{35}.

Suppose for a moment, that the $K$-cycles (5.5) were well defined for all $t\geq 0$. 
Then for $t=0$ the $K$-cycle $(\ell^2(\Delta^R_*),\,\partial\,+\,h)$ would represent $1\in KK_\Gamma(\mathbb{C},\mathbb{C})$ (because $\partial$ is strictly equivariant), and the continuous family (5.5) would provide the desired homotopy between $\gamma$ and the unit $K$-cycle.
\\ 
As we know, this is too much to hope for, because many hyperbolic groups have Kazhdan's property (T), which rules out the existence of such a homotopy. However, according to Lafforgue, hyperbolic groups do not satisfy his strengthened property (T). So one may still hope to find a homotopy as above among Hilbert spaces with $\Gamma$-action of small exponential growth. 
Lafforgue's main theorem states that this is indeed the case:

\begin{theo} (Lafforgue, \cite{41})
Let $(\Gamma,S)$ be a word-hyperbolic group, let $R>0$ be large enough that the Rips complex $\Delta^R_*(\Gamma,d_S)$ is contractible, and let $x_0\in\Delta_0$ be a base point. Fix $\lambda>1$. Then there exist
\begin{itemize}
\item a Hilbert space $\mathcal{H}_{x_0,\lambda}\,=\,\overline{{\mathbb C}(\Delta_*^R(\Gamma,d_S))}$,
given by a completion of the Rips chain complex, \\
\item a hyperbolic distance $d'$ on $\Gamma$ such that $d'-d$ is bounded, \\
\item a contracting square zero chain homotopy of the Rips complex, i.e. a linear map 
$$
h_{x_0}:{\mathbb C}(\Delta_*^R(\Gamma,d_S))\to {\mathbb C}(\Delta_{*+1}^R(\Gamma,d_S))
$$ 
satisfying 
$$
\begin{array}{ccc}
h_{x_0}^2=0, & \partial\circ h_{x_0}\,+\,h_{x_0}\circ\partial\,=\,Id\,-\,p_{x_0}, &
 Im(p_{x_0})\,=\,{\mathbb C}x_0,
\end{array}
\eqno(5.6)
$$
\end{itemize}
such that the following holds:
\begin{itemize}
\item[a)]The maps $F_t\,=\,e^{t d'_{x_0}} (\partial\,+\,h_{x_0})  e^{-t d'_{x_0}}$, where $d'_{x_0}:\Delta^R_*(\Gamma,d_S)\to\mathbb{R}_+$ denotes the distance from the base point, extend to a continuous family of bounded linear operators on $\mathcal{H}_{x_0,\lambda}$.
\item[b)] The natural action of $\Gamma$ on ${\mathbb C}(\Delta_*^R(\Gamma,d_S))$ extends to a continuous representation $\pi$ of exponent $\lambda$ on $\mathcal{H}_{x_0,\lambda}$.
\item[c)] The operators $[F_t,\pi(g)]$ are compact for all $g\in\Gamma$ and all $t\in\mathbb{R}_+$.
\end{itemize}
In particular, the generalized $K$-cycles
$$
{\mathcal E}_t\,=\,\left(\mathcal{H}_{x_0,\lambda},\,e^{t d'_{x_0}} (\partial\,+\,h_{x_0})  e^{-t d'_{x_0}}\right)
\eqno(5.7)
$$
define an exponent-$\lambda$-homotopy between $1\in KK_\Gamma(\mathbb{C},\mathbb{C})$ and $\gamma\in KK_\Gamma(\mathbb{C},\mathbb{C})$.
\end{theo}

The key point is the existence of the desired homotopy for {\bf all} $\lambda>1\,!$
Now one has left the framework of Kasparov's bivariant $K$-theory, but an argument of Higson (\cite{40},(2.12)) shows that the previous theorem still implies ${\bf BC_{coeff}}$. Thus 

\begin{coro} (Lafforgue, \cite{41})
The Baum-Connes conjecture with coefficients is true for all word-hyperbolic groups.
\end{coro}

The demonstration of Lafforgue's theorem requires almost 200 pages and is extremely complicated. We therefore can only outline the strategy of the proof and have to refer to the original paper \cite{41} for details.

\subsection{The case of free groups}

We will study free groups first because the proof of Lafforgue's theorem for free groups is easy and suggests the right strategy for the general case. The fact that $\gamma\,=\,1$ for free groups is due to P.~Julg and A.~Valette \cite{29}. Their work inspired the line of thought followed here.
\\
Let $\Gamma=F_2$ be the nonabelian free group on two generators $s,t,$ and let $S=\{s^{\pm 1},t^{\pm 1}\}$. The geometric realization of the Rips complex $X_*=\Delta^R_*,\,R=1,$ is a tree. Once a base point has been chosen, for example $x_0=x_e\in X_0=F_2,$ there is a canonical contracting simplicial homotopy $h_{x_0}$ of the tree $X_*$: every vertex is sent to the unique simplicial geodesic joining it to the origin. The operators $[F_t,\pi(g)]$ for this homotopy and for the original metric $d=d'$ are of finite rank and will be compact, once they are bounded. The whole problem therefore boils down to find the right Hilbert space.
\\
Rewrite the contracting homotopy as $h_{x_0}=\underset{r}{\sum}\,h_{x_0,r},$ where $h_{x_0,r}:\mathbb{C}(X_0)\to \mathbb{C}(X_1)$ sends a vertex to the edge at distance $r$ on its geodesic journey to the origin $x_0$. By definition $h_{x_0,r}(x_g)=0$ if $\ell(g)<r$. The first step is to replace the $\ell^2$-norm on $\mathbb{C}(X_0)$ by the graph norms of the operators 
$h_{x_0,r}:\ell^2(X_0)\to\ell^2(X_1)$. For $\lambda>1$ and $f\in\mathbb{C}(X_0)$ one puts
$$
\parallel f\parallel_{x_0,\lambda,prel}^2\,=\,\parallel f\parallel_{\ell^2(X_0)}^2\,+\,
\underset{r=1}{\overset{\infty}{\sum}}\,\lambda^{2r}\parallel h_{x_0,r}(f)\parallel_{\ell^2(X_1)}^2,
\eqno(5.8)
$$
(note that the sum is finite) and gains the boundedness of $h:\mathcal{H}_{x,\lambda,prel}\to\ell^2(X_1)$. Lafforgue gives a geometric description of a closely related Hilbert space, which applies immediately to general hyperbolic groups.
Let $e_x,\,x\in X_*,$ be the canonical basis of $\mathbb{C}(X_*)$ and let $l_x,\,x\in X_*,$ be the dual basis. The operator $h_{x,1}$ provides an identification $e:X_0-\{x_0\}\overset{\simeq}{\to} X_1$.
The norm (5.8) can then be rewritten as 
$$
\parallel - \parallel_{x_0,\lambda,prel}^2\,=\,\underset{x\in X_0}{\sum}\vert l_x\vert^2\,+\,\underset{y\in X_0-\{x_0\}}{\sum}\,\lambda^{2r}\vert h_{x_0,r}^t(l_{e(y)})\vert^2,
\eqno(5.9)
$$
where $h_{x_0,r}^t$ denotes the transpose of $h_{x_0,r}$. One has
$$
h_{x_0,r}^t(l_{e(y)})\,=\,\underset{v\in Fl_{x_0}(y,r)}{\sum}\,l_v
\eqno(5.10)
$$
where $Fl_{x_0}(y,r)$ is the {\bf flower} based at $y$ of height $r$, i.e. the set of vertices in $X_0$, which lie at distance $r$ from $y$ and pass through $y$ on their journey to the origin. 

An alternative way to describe flowers is the following. Let $B(x_0,k)$ be the ball around $x_0$ of radius $k=d(x_0,y)$. Then every geodesic path from elements $v,v'\in Fl_{x_0}(y,r)$
to a vertex $w$ in $B(x_0,k)$ will pass through $z$. Consequently 
$$
d(v,w)=d(v,y)+d(y,w)=d(v',y)+d(y,w)=d(v',w).
\eqno(5.11)
$$ 
The {\bf flowers} over $B(x_0,k)$,  i.e. the flowers of arbitrary height and based at points 
of the $k$-sphere around $x_0$ appear thus as the equivalence classes of points of $X_0-B(x_0,k)$ with respect to the equivalence relation
$$
\mathcal{R}_k:\,\,\,x\underset{k}{\sim}x'\,\Leftrightarrow\,d(x,z)=d(x',z),\,\,\,\forall z\in B(x_0,k)
\eqno(5.12)
$$
on $X_0$. Denote by $\overline{Y}_{x_0}^{0,k}$ the set of equivalence classes of $\mathcal{R}_k$. These are the flowers over $B(x_0,k)$ and the points of $B(x_0,k)$. The height $r$ of a flower $Z\in\overline{Y}_{x_0}^{0,k}$ equals $r=d(x_0,Z)-k$. Now Lafforgue defines the Hilbert space $\mathcal{H}_{x_0,\lambda,0}$ as the completion of $\mathbb{C}(X_0)$ with respect to the norm
$$
\parallel - \parallel_{x_0,\lambda,0}^2\,=\,\underset{k=0}{\overset{\infty}{\sum}}\,\,
\underset{Z\in\overline{Y}_{x_0}^{0,k}}{\sum}\,\lambda^{2(d(x_0,Z)-k)}\vert\underset{v\in Z}{\sum}\,l_v\vert^2.
\eqno(5.13)
$$
In this formula the sum over the terms satisfying $d(x_0,Z)-k>0$ gives exactly the second term on the right hand side of (5.9) by (5.10), whereas the sum over the other terms equals a constant multiple of the first term of the right hand side of (5.9). In particular, the norms (5.9) and (5.13) are equivalent. \\
Let us have a closer look at the group action on $\mathcal{H}_{x_0,\lambda,0}$. The norm on this Hilbert space is defined purely in terms of the geometry of the Cayley graph (tree) of $(F_2,S)$, but depends heavily on the choice of the base point $x_0$. Calculating the norm of the operator $\pi(g)$ amounts therefore to bound the norm $\parallel -\parallel_{x_0',\lambda,0}$ with respect to the new base point $x_0'=g^{-1}x_0$ in terms of the original norm $\parallel -\parallel_{x_0,\lambda,0}$. To this end one has to express each flower $Z'\in\overline{Y}_{x_0'}^{0,k}$ over a ball around $x_0'$ as a disjoint union of flowers  $Z_i\in\overline{Y}_{x_0}^{0,k_i}$ over balls around $x_0$. Such a decomposition is not unique, and one is interested in decompositions with as few flowers as possible. Let $Z'=Fl_{x_0'}(y,r)$ be a flower based at $y$. If $y$ does not lie on the geodesic segment $geod(x_0,x_0')$ joining $x_0$ and $x_0'$, then $Z'=Fl_{x_0'}(y,r)=Fl_{x_0}(y,r)$ is simultaneously a flower over balls around $x_0$ and $x_0'$. If $y\in geod(x_0,x_0')$, then 
$$
Z'\,=\,\underset{j}{\coprod}Z_j
\eqno(5.14)
$$
is the disjoint union of at most $C(\Gamma,S)(d(x_0,x_0')+1)$ flowers $Z_j\,=\,Fl_{x_0}(y_j,r_j)$ over balls around $x_0$ whose base point lies at distance 1 from $geod(x_0,x_0')$ and which satisfy 
$$
d(x_0,Z_j)\,\geq\,d(x_0,Z')\,-\,d(x_0,x_0').
\eqno(5.15)
$$ 
Therefore 
$$
\parallel\pi(g)\xi\parallel_{x_0,\lambda,0}^2\,\leq\,C(\Gamma,S)^{\frac12}(1+\ell(g))^{\frac12}\lambda^{\ell(g)}\parallel\xi\parallel_{x_0,\lambda,0}^2,\,\forall\xi\in\mathbb{C}(X),
\eqno(5.16)
$$
i.e. the representation of $F_2$ on $H$ is of exponent $\lambda'$ for every $\lambda'>\lambda$.

\subsection{The Hilbert space, an Ansatz}
The formula (5.13) for the Hilbert space completion of the Rips complex can easily be generalized.
Let $(\Gamma,S)$ be a word-hyperbolic group. Let $R>0$ be large enough that the Rips-complex $\Delta_*=\Delta^R_*(\Gamma,d_S)$ is contractible and fix a base point $x_0\in\Delta_0$. Rips $p$-simplices correspond then to oriented $(p+1)$-element subsets of $\Gamma$ of diameter at most $R$.\\

\begin{figure}[htpb]
\begin{center}
\includegraphics[scale=0.5]{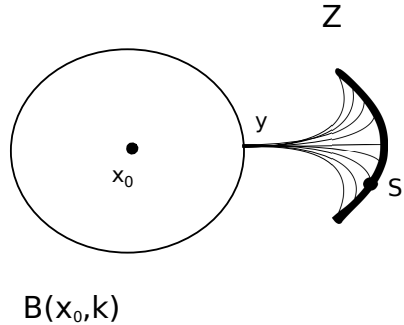}
\caption{A flower in $\overline{Y}_{x_0}^{p,k}$}
\label{default}
\end{center}
\end{figure}

Lafforgue says that two Rips $p$-simplices $S',S''$ are $k$-equivalent, if there exists an isometry between $B(x,k)\cup S'$ and $B(x,k)\cup S''$, which sends $S'$ to $S''$ while preserving orientations, and fixes $B(x,k)$ pointwise. This is an equivalence relation and the set of equivalence classes is denoted by $\overline{Y}_{x}^{p,k}$. The equivalence classes are called {\bf flowers} over $B(x_0,k)$ if $d(x,S')>k$ and equal a single simplex if $d(x,S')\leq k-R$.\\
Lafforgue defines now, similar to (5.13), a Hilbert space $\mathcal{H}_{x_0,\lambda,0}(\Delta_p)$ as the completion of $\mathbb{C}(\Delta_p)$ with respect to the norm 
$$
\parallel -\parallel_{x_0,\lambda,0}^2\,=\,\underset{k=0}{\overset{\infty}{\sum}}\,\,
\underset{Z\in\overline{Y}_{x_0}^{p,k}}{\sum}\,\lambda^{2(d(x_0,Z)-k)}\vert\underset{S'\in Z}{\sum}\,l_{S'}\vert^2,
\eqno(5.17)
$$
where again $(l_S),S\in\Delta_p,$ is the dual of the canonical basis of $\mathbb{C}(\Delta_p)$.\\
The left translation action of $\Gamma$ on the Rips complex gives rise to a representation $\pi$ on $\mathcal{H}_{x_0,\lambda,0}(\Delta_p)$, which is of exponent $\lambda'$ for every $\lambda'>\lambda$ by essentially the same argument as in the case of free groups. A direct consequence is

\begin{theo}(Lafforgue,\cite{38})
Hyperbolic groups do not have property ${\bf (T)_{Hilb}^{Strong}}$.
\end{theo}

This follows easily from the properties of $\mathcal{H}_{x_0,\lambda,0}(\Delta_0)$. It is clear that the regular representation $\pi$ of $\Gamma$ on $\mathcal{H}_{x,\lambda,0}(\Delta_0)$ has no fixed vectors. However, there is a fixed vector for its contragredient representation $\check{\pi}$, because the 
linear functional 
$$
l:\sum a_gx_g\mapsto\,\sum a_g
\eqno(5.18)
$$ 
is bounded: one has 
$$
l\,=\,\underset{r=0}{\overset{\infty}{\sum}}\,\underset{d(x_0,v)=r}{\sum}\,l_v\,=\,\underset{Z\in \overline{Y}_{x_0}^{0,0}}{\sum}\,\underset{v\in Z}{\sum}\,\,l_v,
\eqno(5.19)
$$
(note that the flowers in $\overline{Y}_{x_0}^{0,0}$ are just the spheres around $x_0$), so that 
$$
\vert l(\xi)\vert^2\,=\,\vert
\underset{Z\in\overline{Y}_{x_0}^{0,0}}{\sum}\,\langle\lambda^{-d(x_0,Z)},\,\lambda^{d(x_0,Z)}\,(
 \underset{v\in Z}{\sum}\,l_v)\rangle\vert^2\,\leq\,(1-\lambda^{-2})^{-1}
\parallel \xi\parallel_{x,\lambda,0}^2
\eqno(5.20)
$$
by the Cauchy-Schwarz inequality. Suppose that $\Gamma$ possesses property ${\bf (T)_{Hilb}^{Strong}}$. Then, for $\lambda$ sufficiently close to one, there exists a self-adjoint ``Kazhdan"-projection $p\in \mathcal{C}_\lambda(\Gamma,\ell_S,\mathcal{H})$. It satisfies $\pi(p)=0$, but $\check{\pi}(p)\neq 0$ cannot be zero, because its image contains $l$. This is impossible because $\check{\pi}(p)\,=\,\pi(p)^*$ for self-adjoint projections.

\subsection{Metrically controlled operators}
Lafforgue's Hilbert spaces have interesting properties: they are defined in terms of the geometry of the Cayley-graph, the regular representation on them may be of arbitrary small exponent and it cannot be separated from the trivial representation. This makes them into excellent candidates for the Hilbert spaces needed in Theorem 5.2. \\
However, for general hyperbolic groups, it is difficult to establish the boundedness of 
any contracting chain homotopy of the Rips complex. 
The naive solution of completing the natural domain $\mathbb{C}(\Delta_*)$ of such an operator with respect to the graph norm will not work: contrary to the case of free groups, it will destroy the purely geometric nature of the Hilbert space and one will loose the control over the norm of the representation on it.\\ In order to solve the problem Lafforgue proceeds in two steps:
He identifies a class of operators, which are sufficiently controlled by the geometry of $(\Gamma,d_S)$, so that their graph norms are equivalent to norms of the type (5.17) considered before. A weighted sum of iterated graph norms can then still be defined in a purely geometric way and leads to Hilbert spaces, on which all controlled operators act boundedly. 
Then Lafforgue constructs very carefully controlled contractions of the Rips complex.\\
His notion of control, inspired by the properties of the operators $h_r$ in the case of free groups, is 
\begin{defi} Let $M,r_1,r_2>1$, and let $r\in\mathbb{N}$. 
A linear operator 
$$
\Phi:\mathbb{C}(\Delta_*)\to\mathbb{C}(\Delta_*)
$$ 
is $r_1$-geodesic (with respect to $x_0$) and $(M,r_2)$-controlled of propagation $r$, if its matrix coefficients $\Phi_{S_0,S_1},S_0,S_1\in\Delta_*,$ satisfy the following conditions:
\begin{itemize}
\item $\Phi_{S_0,S_1}=0$ unless $d(S_{1},geod(x_0,S_0))<r_1$ and $\vert d(S_0,S_1)-r\vert<r_2$.
\item If there exists an isometry between 
$$
\begin{array}{ccc}
B(x_0,M)\cup B(S_0,M)\cup B(S_1,M) & \text{and} & B(x_0,M)\cup B(S'_0,M)\cup B(S'_1,M),
\end{array}
$$ 
which sends $S_0$ to $S'_0$, $S_1$ to $S_1'$, while preserving orientations, and fixes $B(x_0,M)$ pointwise, then the matrix coefficients $\Phi_{S_0,S_1}$ and $\Phi_{S'_0,S'_1}$ coincide. 
\item The set of all matrix coefficients is bounded.
\end{itemize}
\end{defi}
This notion suggests the following modification of Lafforgue's Hilbert space.
\begin{defi} 
Fix $M,r_1>1$ and let $k\in\mathbb{N}$. Let $\overline{Y}_{x_0}^{p,k,m}$ be the set of $m$-fold iterated, 
$M$-thickened {\bf flowers} of $p$-simplices over $B(x_0,k)$, i.e. the set of equivalence classes of $(m+1)$-tuples $(S_0,S_1,\ldots,S_m),m\in\mathbb{N},$ of Rips-simplices $S_i\in\Delta_*,\,S_0\in\Delta_p,$ such that $d(S_{i+1},geod(x_0,S_i))<r_1$ and $d(x_0,S_m)>k+2M$, with respect to the following equivalence relation. Two $(m+1)$-tuples 
$(S_0,\ldots,S_m)$ and $(S_0',\ldots,S_m')$ are equivalent if there exists an isometry between the subsets 
$$
B(x_0,k+2M)\cup B(S_0,M)\cup\ldots\cup B(S_m,M)
$$ 
and 
$$
B(x_0,k+2M)\cup B(S'_0,M)\cup\ldots\cup B(S'_m,M)
$$ 
of $(\Gamma,d_S)$, which maps $S_i$ to $S_i'$ for all $i$ (preserving the orientations of $S_0$ and $S'_0$) and fixes $B(x_0,k+2M)$ pointwise. 
\end{defi}

\begin{figure}[htpb]
\begin{center}
\includegraphics[scale=0.5]{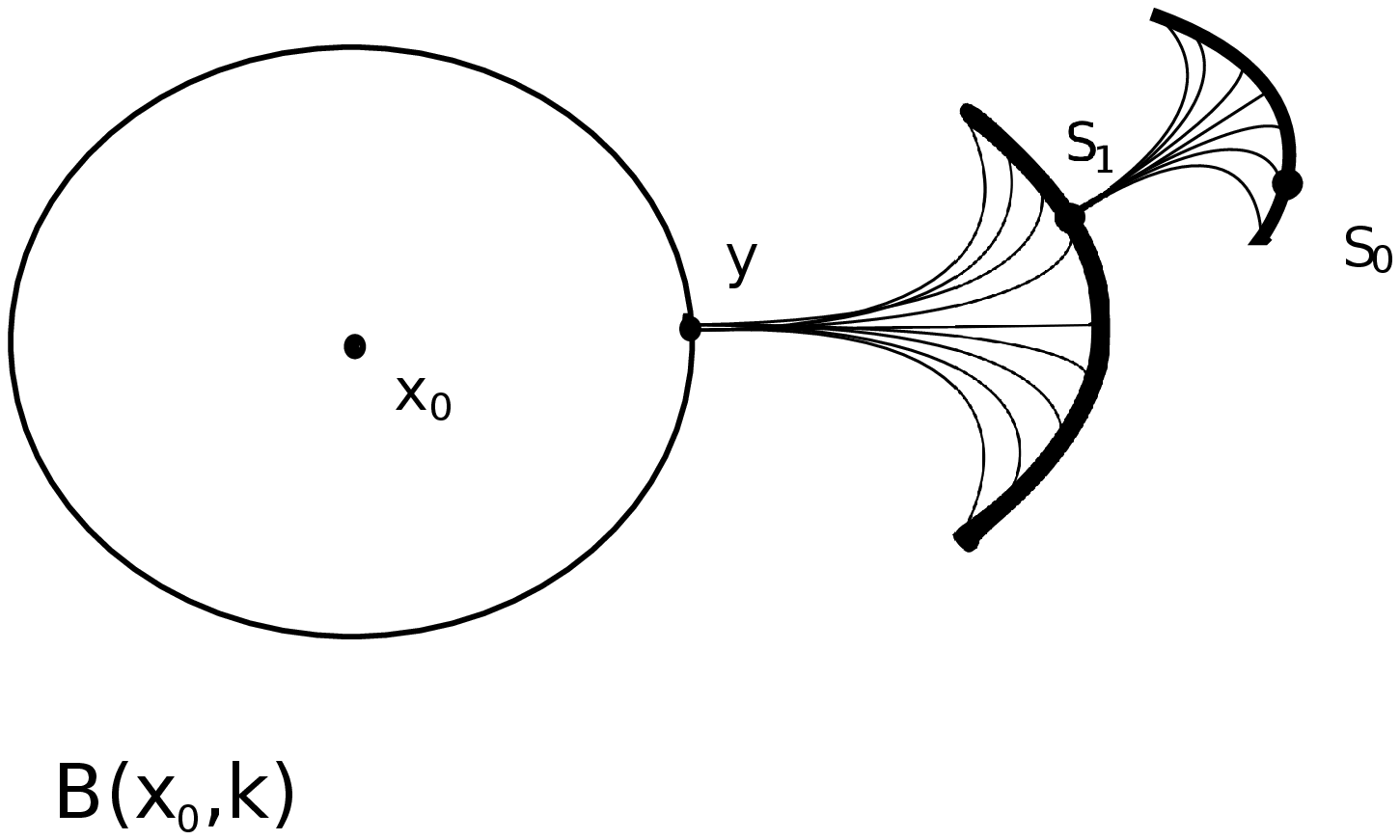}
\caption{A flower in $\overline{Y}_{x_0}^{p,k,1}$}
\label{default}
\end{center}
\end{figure}

Lafforgue defines then the Hilbert space $\mathcal{H}_{x_0,\lambda}(\Delta_p)$ as the completion of $\mathbb{C}(\Delta_p)$ with respect to the norm
$$
\parallel-\parallel_{x_0,\lambda}^2\,=\,\underset{m=0}{\overset{\infty}{\sum}}\,B^{-m}\left(\underset{k=0}{\overset{\infty}{\sum}}\,\,
\underset{Z\in\overline{Y}_{x_0}^{p,k,m}}{\sum}\,\lambda^{2(d(x_0,S_0)-k)}\vert\underset{(S_0,\ldots,S_m)\in Z}{\sum}\,l_{S_0}\vert^2\right),
\eqno(5.21)
$$
where $B=B(\lambda)>>0$ is large, but fixed. Its new basic property is described in

\begin{lemm}
Every linear map $\Phi:\mathbb{C}(\Delta^R_p)\to\mathbb{C}(\Delta^R_q)$, which is $r_1$-geodesic (with respect to $x_0$) and metrically $(M,r_2)$-controlled and of fixed bounded propagation, extends to a bounded linear operator 
$$
\Phi:\mathcal{H}_{x_0,\lambda}(\Delta^R_p)\to \mathcal{H}_{x_0,\lambda}(\Delta^R_q).
$$
\end{lemm}
In fact, let $Z_1\in\overline{Y}_{x_0}^{p,k,m}$. Then, as $\Phi$ is $r_1$-geodesic and $(M,r_2)$-controlled of fixed propagation, its transpose $\Phi^t$ satisfies
$$
\Phi^t(\underset{(S_1,\ldots,S_m)\in Z_1}{\sum}\,l_{S_1})\,=\,\underset{Z}{\sum}\,\,\alpha_{Z,Z_1}
\left( \underset{(S_0,S_1,\ldots,S_m)\in Z}{\sum}\,l_{S_0}\right),
\eqno(5.22)
$$
where
\begin{itemize}
\item The $(m+1)$-fold iterated $M$-thickened flowers $Z$ ``prolongate" the $m$-fold iterated $M$-thickened flower $Z_1$.
\item The number of flowers $Z$, which occur in the sum on the right hand side, is bounded by an absolute constant $C_1(\Gamma,S,R,M,r_1,r_2)$. 
\item $\alpha_{Z,Z_1}\,=\,\Phi_{S_0,S_1}$ is a matrix coefficient of $\Phi$, which depends only on $Z$.
\end{itemize}
Thus 
$$
\parallel\Phi(\xi)\parallel_{x_0,\lambda}\,\leq\,C_2\,\parallel\xi\parallel_{x_0,\lambda}
\eqno(5.23)
$$
by the Cauchy-Schwarz inequality, where $C_2$ depends on $C_1$, the $\ell^\infty$-norm of the matrix coefficients of $\Phi$, and of the propagation $r$ of $\Phi$. \\
\\
With this result at hand, the next step is to look whether the homotopy (5.7) may be realized by geodesic and metrically controlled operators. The simplicial differential of the Rips complex is obviously $R$-geodesic, $R$-controlled, and of propagation at most $R$. Concerning the contracting chain homotopy $h_{x_0}$, the story is more complicated.
\\
There is a standard procedure for contracting the Rips complex of a hyperbolic group \cite{18},\cite{33},\cite{41}. For a given Rips $p$-simplex $S_0$, one lets $\widetilde{h}_{x_0}'(S_0)$ be a mean over the Rips $(p+1)$-simplices $S_0\cup\{y\}$ with $y\in\Gamma$ of minimal word-length (distance to the origin). Put $\psi_{x_0}=Id-(\widetilde{h}_{x_0}\partial+\partial \widetilde{h}_{x_0})$. Then $\underset{n=0}{\overset{\infty}{\sum}}\,\widetilde{h}_{x_0}\psi_{x_0}^n$ will be a contracting chain-homotopy of $\mathbb{C}(\Delta_*)$. The construction shows also that the subcomplexes spanned by Rips simplices supported in a given ball $B\subset(\Gamma,d_S)$, are contractible as well (take the center of the ball as new origin).
\\
The operators $h_{x_0}'$ and $\psi_{x_0}$ are clearly geodesic, controlled, and of finite propagation. So they extend to bounded operators on $\mathcal{H}_{x_0,\lambda}(\Delta_*)$. However, one cannot prevent the norms of the powers $\psi_{x_0}^n,\,n\in\mathbb{N},$ to grow exponentially ! This is due to the fact, that for fixed control parameters the operator $\psi_{x_0}^n$ will not anymore be geodesic and controlled if $n$ is large: the transition coefficient from a simplex $S_0$ to $S_1$ will depend on the whole trajectory from $S_0$ to $S_1$, and not only on the geometry of the union of  balls of fixed radius around $x$, $S_0$ and $S_1$.
\\
Lafforgue solves the problem by constructing ad hoc geodesic and metrically controlled chain maps 
$$
\varphi_{x_0,r}:\mathbb{C}(\Delta_*)\to \mathbb{C}(\Delta_*),
\eqno(5.24)
$$ 
which cover the identity in degree -1 and move each simplex $r$ steps towards the origin.\\
\\
For $\{y\}\in\Delta_0$ let $\varphi_{x_0,r}(\{y\})$ be the mean over the points of $geod(x_0,y)\cap S(y,r),$ and extend by linearity to $\mathbb{C}(\Delta_0)$. If $\varphi_{x_0,r}$ has been defined on $\mathbb{C}(\Delta_k)$, and if $S'$ is a $(k+1)$-Rips-simplex, then let $\varphi_{x_0,r}(S')$ be a filling of the cycle $\varphi_{x_0,r}(\partial S')$ inside a fixed ball of diameter $R+2\delta$ containing $\underset{y\in S'}{\bigcup}\,geod(x_0,S')\cap S(y,r),$ (such  fillings exist as remarked before), and extend by linearity to $\mathbb{C}(\Delta_k)$. Finally use the same 
procedure to construct geodesic and metrically controlled homotopy operators
$$
h'_{x_0,r}:\mathbb{C}(\Delta_*)\to \mathbb{C}(\Delta_{*+1})
\eqno(5.25)
$$ 
satisfying 
$$
\varphi_{x_0,r+1}-\varphi_{x_0,r}\,=\,h'_{x_0,r}\partial+\partial h'_{x_0,r},
\eqno(5.26)
$$ 
and put 
$$
h'_{x_0}=\underset{r}{\sum}\,h'_{x_0,r}
\eqno(5.27)
$$
Lafforgue combines the two constructions to obtain a contracting square zero homotopy 
$$
h''_{x_0}\,=\,h'_{x_0}\circ\psi_{x_0}^N\,+\,\underset{n=0}{\overset{N-1}{\sum}}\,\widetilde{h}_{x_0}\circ\psi_{x_0}^n,
\eqno(5.28)
$$
which is controlled and moves simplices strictly towards the origin. The linear maps
$$
F_t\,=\,e^{t d_{x_0}} (\partial\,+\,h''_{x_0})  e^{-t d_{x_0}}
\eqno(5.29)
$$
extend then to bounded operators on $\mathcal{H}_{x_0,\lambda}(\Delta_*)$
for every $\lambda>1$ and every $t\in\mathbb{R}_+$.
 
\subsection{The pigeonhole principle and the end of the proof}
What remains to be done ? Still one step, and it is by far the hardest: it has to be shown that 
the commutators $[F_t,\pi(g)]$ are compact for all $g\in\Gamma$ and $t\in\mathbb{R}_+$. 
The work of Kasparov-Skandalis \cite{33} and of Mineyev-Yu \cite{46} suggests, how to proceed:

\begin{itemize}
\item The contracting homotopy $h''_{x_0}$ has to be replaced by a contraction $h_{x_0}$ which depends ``continuously" on $x_0$:
$$
\underset{S_0\to\infty}{\lim}\,\left(h_{x_0}(S_0)-h_{x'_0}(S_0)\right)\,=\,0,\,\forall x'_0\in\Delta_0.
\eqno(5.30)
$$
\item The word-metric $d$ on $\Gamma$ has to be changed into an equivalent ``continuous metric" $d'$ with the crucial property 
$$
\underset{x\to\infty}{\lim}\,\left(d'(x_0,x)\,-\,d'(x_0',x)\right)\,=\,0,\,\forall x'_0\in\Delta_0.
\eqno(5.31)
$$ 
\end{itemize}

The necessary changes are quite subtle. They bring a pigeonhole argument into play which allows to deduce (5.30) and (5.31). Unfortunately, the operator $h_{x_0}$ and the operator 
$d'_x:\mathbb{C}(\Delta^R_*)\to \mathbb{C}(\Delta^R_*)$ of multiplication with the distance from the origin cannot be metrically controlled anymore. \\

This forces Lafforgue to modify again the underlying Hilbert space. The final formula for the norm turns out to be much more complicated than our ``baby model" (5.21) and depends on nine parameters (instead of the four parameters $R,M,B,r_1$ we used). We will only give some indications and refer to Lafforgue's original paper for more detailed information.
\\
In the notations of 5.2 it suffices more or less to verify that the operators $[\pi(g),h_{x_0,r}]$ and $[\pi(g),d'_{x_0}]$
are compact for every $r\geq 0$ and $g\in S$ (and thus for every $g\in G$). (The commutators $[\pi(g),\partial]$ vanish because the differential $\partial$ is $\Gamma$-equivariant.) In more convenient terms, this means that the operators
$$
h_{x_0',r}-h_{x_0,r},\,\,\,\text{and}\,\,\,d'_{x_0}-d'_{x'_0},\,\,\, x_0,x_0'\in\Delta_0,\,\,\, d(x_0,x_0')=1,
\eqno(5.32)
$$ 
are compact.\\
\\
Decompose the operator (5.28) into a sum of operators of propagation $r$:  $h_{x_0}''=\underset{r}{\sum}h_{x_0,r}''$. The construction of the operator $h''_{x_0,r}$ depends on a large number of choices. Any of these choices was sufficient to arrive at (5.29), but now, one has good reason to keep track of the choices made. So Lafforgue introduces a probability space $(\Omega,\mu)$, whose points label the possible choices in the construction of $h''_{x_0,r},r\in{\mathbb N}$.
The corresponding operators are denoted by $h''_{x_0,r,\omega},\omega\in\Omega$.
If a simplex $S_0\in\Delta_*$ is very far from the origins $x_0$ and $x_0'$ in the sense that 
$d(x_0,S)>>r$, then one may hope that $h''_{x_0,r,\omega}(S_0)=h''_{x_0',r,\omega}(S_0)$ with a high probability. In fact
\begin{lemm} 
There exists a universal constant $C>0$ such that 
$$
\mu(\{\,\omega\in\Omega,\,h''_{x_0,r,\omega}(S_0)\,\neq\,h''_{x_0',r,\omega}(S_0)\,\})\,\leq\,
\frac{C\,d(x_0,x_0')}{1+(d(x_0,S_0)-r)}
\eqno(5.33)
$$
for all $S_0\in\Delta_*$, all $r\in\mathbb{N}$, and all $x_0,x_0'\in\Delta_0=\Gamma$.
\end{lemm}
The proof uses a counting argument as in \cite{33}, which is based on a pigeonhole principle. The operator
$$
h_{x_0}\,=\,\underset{r=0}{\sum}\,h_{x_0,r},\,\,\,h_{x_0,r}\,=\,\underset{\Omega}{\int}\,h''_{x_0,r,\omega}\,d\mu,
\eqno(5.34)
$$
will be the definitive contracting homotopy required in 5.2. Suppose that the operators $h''_{x,r,\omega},\omega\in\Omega,$ were uniformly bounded with norms summable with respect to $r$.  
The operator $h_{x}-h_{x'}$ would then split as a sum of an operator of arbitrary small norm (neglecting simplices close to the origin) and a finite rank operator. In other words, it would be compact.\\
In the same spirit, Lafforgue introduces a family of modified metrics $d_{\widetilde{\omega}},\,\widetilde{\omega}\in\widetilde{\Omega},$ labeled by another probability space $(\widetilde{\Omega},\widetilde{\mu})$. They are obtained by an extremely subtle averaging process, and satisfy
$$
\widetilde{\mu}(\{\,\widetilde{\omega}\in\widetilde{\Omega},\,d_{\widetilde{\omega}}(x_0,x)\,\neq\,d_{\widetilde{\omega}}(x_0',x)\,\})\,\leq\,
\frac{C(d(x_0,x_0'))}{1+d(x_0,x)}.
\eqno(5.35)
$$
Consequently the final new metric 
$$
d'\,=\,\underset{\widetilde{\Omega}}{\int}\,d_{\widetilde{\omega}}\,d\widetilde{\mu},
\eqno(5.36)
$$
used in 5.2, is ``continuous" over large distances, and one could essentially conclude 5.2.c), provided that the operators $d_{x_0,\omega'}-d_{x_0',\omega'}$ were uniformly bounded.\\
However, one cannot prove this: neither the operators $h''_{x,r,\omega},$ nor the operators $d_{x_0}-d_{x_0,\widetilde{\omega}}$ are metrically controlled! \\
\\
Fortunately, they almost are: their matrix coefficients $a_{S_0,S_1}$ vanish unless\\ $d(S_1,geod(x,S_0))<r_1$ and $\vert d(S_0,S_1)-r\vert<r_2$ (resp. $S_0=S_1$), and depend only on the isometry class of the union of $B(x,M)\cup B(S_0,M)\cup B(S_1,M)$ and a uniformly finite family of ``control sets" of uniformly bounded diameter, located along $geod(x,S_0)$.\\
This leads Lafforgue to the definitive version of his Hilbert space. He enriches the definition of the iterated flowers 5.6 by the introduction of uniformly finite families of ``control sets" of uniformly bounded diameter, located uniformly close to $geod(x,S_0)$. The corresponding isometry relation has to take control sets into account and the number of control sets introduces a weight factor in the sums defining the Hilbert norm. Another weight factor, given by a very mildly decaying exponential of the cardinality of each flower has still to be introduced to arrive finally at a Hilbert space satisfying (5.7). The proof of theorem 5.2 is thus complete.

\subsection{Concluding remarks}

Lafforgue's proof of the Baum-Connes conjecture with coefficients for word-hyperbolic groups shows once more the power and flexibility of Kasparov's ``Dirac-Dual Dirac" approach. It may even work in the presence of Kazhdan's property (T) as he shows. His work on the strengthened property (T) indicates however, that the method has its limits. There seems to be no way to apply it in the crucial case of lattices in simple algebraic groups of split rank $\geq 2$ over local fields. At present the search for a ``truly noncommutative" version of J.~B.~Bost's Oka-principle \cite{10} seems to be the only hope to settle this case. New ideas will be needed to decide whether the fascinating predictions of Baum and Connes hold for further classes of groups.

\end{document}